\documentclass[3p,times]{elsarticle}
\usepackage[utf8]{inputenc}
\usepackage[english]{babel}
\usepackage{amssymb,amsmath}
\usepackage{gensymb}
\usepackage{stmaryrd}
\usepackage{graphicx}
\usepackage{float}
\usepackage{algorithm}
\usepackage{algpseudocode}
\usepackage{cclicenses}
\usepackage[textsize=footnotesize]{todonotes}

\def\be{\boldsymbol e}
\def\bx{\boldsymbol x}

\def\calG{\mathcal G}
\def\calM{\mathcal M}
\def\polI{\mathbb I}
\def\calD{\mathcal D}
\def\calL{\mathcal L}

\newcommand{\ie}{{\em i.\thinspace{}e. }}

\newcommand{\eg}{{\em e.\thinspace{}g. }}

\newcommand{\ds}{\displaystyle}
\newcommand{\pl}{\partial }

\newcommand{\Ste}{\textit{Ste}}

\newcommand{\eps}{\varepsilon}

\newcommand{\DIV}[1][]{\nabla_{#1}\cdot}
\def\divt{\DIV[\tilde x]}
\newcommand{\eigrad}[1][\be_i]{#1\!\cdot\!\nabla_{\tilde x}}

\usepackage{color}

\DeclareMathOperator\erf{erf}
\DeclareMathOperator\erfc{erfc}

\journal{International Journal of Heat and Mass Transfer}
 \date{March 27, 2025}

\begin{document}


\begin{frontmatter}

\title{An implicit regularized enthalpy Lattice Boltzmann Method for the Stefan problem}

\author[lmrs]{Francky Luddens\corref{ca}}
\ead{francky.luddens@univ-rouen.fr}

\author[larema]{Corentin  Lothod{\'e}}
\ead{corentin.lothode@univ-angers.fr}

\author[lmrs]{Ionut Danaila}
\ead{ionut.danaila@univ-rouen.fr}

\address[lmrs]{Univ Rouen Normandie, CNRS, Normandie Univ, LMRS UMR 6085, F-76000 Rouen, France}
\address[larema]{Univ Angers, CNRS, LAREMA - UMR 6093, SFR MATHSTIC, F-49000 Angers, France}

\cortext[ca]{Corresponding author.}


\begin{abstract}
Solving the Stefan problem, also referred as the heat conduction problem with phase change, is a necessary step to solve phase change problems with convection.
In this article, we are interested in using the Lattice Boltzmann Method (LBM) to solve the Stefan problem using a regularized total enthalpy model. The liquid fraction is treated as a non-linear source/sink term, that involves the time derivative of the solution. The resulting non-linear system is solved using a Newton algorithm.
By conserving the locality of the problem, this method is highly scalable, while keeping a high accuracy.
The newly developed scheme is analyzed theoretically through a Chapman-Enskog expansion and illustrated numerically with 1D and 2D benchmarks.
\end{abstract}

\end{frontmatter}
\textit{License:} This work is licensed under \texttt{https://creativecommons.org/licenses/by/4.0/} \includegraphics[height=1em]{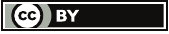}

\section{Introduction}


In this article, we consider the classical two-phase Stefan problem \cite{lame1831,tarzia2000}.
It describes the melting or solidification of a material taking into account the conduction as the main heat transfer mechanism.
Numerical solvers for liquid-solid phase change systems, such as those occurring in metal casting, energy storage, cryosurgery or food freezing, are usually first benchmarked against the analytical solution of the Stefan problem \cite{kowalewski2004phase}.

The Stefan problem is well known in the mathematical community, with many theoretical results published, see for example \cite{rubinstein1979stefan} for a review.
In the context of numerical simulation, several frameworks were developed to solve this problem \cite{voller1990fixed}.
In the finite element community, different approaches were suggested, from enthalpy based methods (e.g. \cite{elliott1987error}) to front tracking methods (e.g. \cite{li1983finite}).
For other numerical frameworks, similar approaches are developed, see for example \cite{muhieddine2009various} for finite volumes and \cite{savovic2003finite} for finite differences.

The Lattice Boltzmann Method (LBM) was first adapted in~\cite{wolf1995lattice} to tackle heat diffusion.
Later, a chemical reaction analogy and two types of quasi particles were introduced in~\cite{de1998mesoscopic}. The resulting method involves solving a linear system of equations that was subsequently simplified in~\cite{miller2001lattice} using only one type of quasi particle and a phase field method.
The position of the interface is recovered through the liquid fraction.
Another approach is to use a source term and an iterative procedure as proposed initially in~\cite{jiaung2001lattice}.
In their work, they also derive a Chapman-Enskog analysis of the scheme.
It was latter extended to phase change problems with convection, see \cite{chatterjee2005enthalpy,chatterjee2008enthalpy}.
The same model is used in~\cite{huber2008lattice}, with only one iteration for the source term at the expense of a lower accuracy.
An implicit method is derived in~\cite{eshraghi2012implicit}, introducing a range of temperature for melting to model binary alloys which can be seen as a regularization of the enthalpy.
Another approach consists in modifying the equilibrium function to recover enthalpy as proposed by \cite{huang2013new}.
Authors of \cite{luo2015lattice} compare the hybrid LBM finite differences scheme introduced in \cite{lallemand2003hybrid} and the enthalpy approach in the context of phase change materials with convection.
They found that the enthalpy method leads to more accurate results.

As noted in \cite{voller1981accurate}, the enthalpy method generates two non-physical effects when dealing with pure materials with a single fusion/solidification temperature: (i) the phase change interface increases step by step like a staircase; (ii) the temperature history at a given spatial point leads to a number of rogue plateaux in addition to the physical temperature plateau.

Following the recent work of \cite{sadaka2020parallel}, we propose another formulation of the source term, mixing a local Newton algorithm and a regularized enthalpy method.
This approach has already shown promising results in the context of numerical methods based on finite elements.
Adapting this method to LBM while keeping efficient parallelization is not a trivial task since it involves solving a non-linear problem during the stream-collide step.

The paper is organized as follows: in \S\ref{sec:problem} we describe the Stefan problem. \S\ref{sec:mesh} provides a description of the discretization, and the link between a mesh (e.g. suitable for finite differences schemes) and the lattice is expressed. We also present previous methods to handle the enthalpy term.
The novelty of this paper is described in \S\ref{sec:Newton_enthalpy}, with a Newton method to compute the regularized enthalpy as a source term.
\S\ref{sec:analysis} is devoted to the numerical analysis of the method, using the Chapman-Enskog expansion.
Finally, \S\ref{sec:numerical} contains numerical illustrations in 1D and 2D to assess the stability and accuracy of the method.

\section{Stefan problem and total enthalpy model}\label{sec:problem}

Consider the melting of a semi-infinite one-dimensional pure material, initially at uniform temperature $T_0 < T_f$, where $T_f$ is the fusion temperature. At the left border $x=0$, the temperature is suddenly increased to $T_h >T_f$ and maintained constant for any time instant $t>0$. The liquid phase starts to propagate from left to right. Solid and liquid domains are denoted by subscripts $s$ and $l$, respectively. They are separated by the interface located at $x_f(t)$. We use the following scaling:
\begin{itemize}
\item $x \rightarrow x/L$,
\item $t \rightarrow t/(L^2/\alpha_l)$,
\item $\theta =(T-T_f)/(T_h-T_f)$,
\end{itemize}

 \noindent with $L$ a reference length scale and $\alpha_l=k_l/(\rho_l c_l)$ the thermal diffusivity. Density $\rho$, specific heat $c$ and thermal conductivity $k$ are constants, but different from one phase to another.

The Stefan problem is completely defined by writing the corresponding heat equation in each domain, with appropriate boundary conditions, and by imposing the balance of heat fluxes at the interface (Stefan condition). In the non-dimensional setting, the problem formulation is summarized in Eq. \eqref{eq-Stefan}, where $h_{sl}$ is the latent heat and $\Ste$ the Stefan number.
\begin{equation}
	\begin{tabular}{l|l}
		{Liquid phase} & 	{Solid phase} \\ 
		$\ds \frac{\pl \theta_l}{\pl t} =  \frac{\pl^2 \theta_l}{\pl x^2} $& $\ds \frac{\pl \theta_s}{\pl t} =  \alpha \frac{\pl^2 \theta_s}{\pl x^2}, \quad \alpha=\frac{\alpha_s}{\alpha_l},$ \\[0.2cm] 
		$\ds \theta_l(0,t) =1, \quad \theta_l(x_f,t) =0$ & $\ds \theta_s(x_f,t) =0, \quad \theta_s(+\infty,t) =\theta_{0} < 0$ \\[0.2cm] \hline
		\multicolumn{2}{c}{Stefan condition at the interface $x=x_f$} \\[0.2cm]
		\multicolumn{2}{c}{$\ds - \frac{\pl \theta_l}{\pl x} + k  \frac{\pl \theta_s}{\pl x} = \frac{\rho}{\Ste}  \frac{dx_f}{dt},\quad \mbox{with}
			\quad \rho=\frac{\rho_s}{\rho_l}, \quad k=\frac{k_s}{k_l}, \quad \Ste = \frac{c_{l} (T_h-T_f)}{h_{sl}}.$}
	\end{tabular}
	\label{eq-Stefan}
\end{equation}
The following solution to Eq. \eqref{eq-Stefan} is obtained by assuming that the interface is moving following the law $x_f= 2 \lambda \sqrt{t}$, with $\lambda$ a constant:
\begin{equation}\label{eq-exact-sol}
	\theta_{l}(x,t) = 1-\frac{1}{\erf(\lambda)}\erf\left(\frac{x}{2\sqrt{t}}\right), \quad
	\theta_{s}(x,t) = \theta_{0}-\frac{\theta_{0}}{\erfc\left(\ds \frac{\lambda}{\sqrt{\alpha}}\right)} \erfc\left(\ds \frac{x}{2\sqrt{\alpha t}}\right),
\end{equation}
where $\erf(z) =(2/\sqrt{\pi}) \int_0^z \exp(-\tau^2) d\tau$ is the error function and $\erfc(z)=1-\erf(z)$. The constant $\lambda$
is determined from the Stefan condition, \ie by solving the following equation:
\begin{equation}
	\frac{\exp\left(-\lambda^2\right)}{\erf(\lambda)}+ \theta_{0}\frac{k}{\sqrt{\alpha}}\frac{\exp\left(-\frac{\lambda^2}{\alpha}\right)}{\ds \erfc\left(\frac{\lambda}{\sqrt{\alpha}}\right)}=\frac{\rho\lambda\sqrt{\pi}}{\Ste}.
	\label{eq-Stefan-lambda}
\end{equation}

The total enthalpy  model used to solve numerically the Stefan problem \eqref{eq-Stefan} is based on a temperature transformation of the energy conservation law. For the sake of clarity, we only present the case where the materials properties are similar in the liquid phase and the solid phase (hence $\rho = k = \alpha= 1$). The enthalpy $H$ is written as  the sum of the sensible heat and the latent heat:
$$H(\theta) = \theta + \dfrac{\varphi(\theta)}{\Ste} = \theta + \begin{cases} 0 &\textrm{ if } \theta \leq 0, \\ \dfrac{1}{\Ste} &\textrm{ if } \theta > 0,\end{cases}$$
with  $\varphi$ the liquid fraction (taking the value $1$ if $\theta>0$ and $0$ if $\theta \leq 0$). 
The following equation is obtained
\begin{equation}
	\frac{\partial H}{\partial t}  -
	\Delta \theta   = 0,
	\label{eq-H}
\end{equation}
which is equivalent to
\begin{equation}
	\frac{\partial \theta}{\partial t}  -
	\Delta \theta  +  \frac{1}{\Ste}\frac{\partial  \varphi(\theta)}{\partial t}  = 0.
	\label{eq-H-source-sink}
\end{equation}

Note that the formulations~\eqref{eq-H} or~\eqref{eq-H-source-sink} of the Stefan problem are valid in 2D and 3D for phase change problems without convection.

\section{LBM methods for the Stefan problem}\label{sec:mesh}

\subsection{Mesh and lattice}
In this note, we will use the LBM to solve the PDE~\eqref{eq-H-source-sink}.
When needed, we will link the lattice to an actual mesh of the domain $\calD$, resulting in a uniform Cartesian mesh.
We denote by $N$ the size of the lattice, $\delta x = \frac{1}{N}$, and $\delta t$ the time step.
We define the lattice $\calL=\left[0,N\right]^d$, where $d$ is the space dimension.
To properly perform the Chapman-Enskog analysis, we introduce the following rescaling:
$$
t_L := \frac{t}{\delta t},\quad \bx_L:= \frac{\bx}{\delta x},\quad \theta_L(t_L,\bx_L):=\theta(t,\bx),\quad H_L(t_L,\bx_L):=H(t,\bx).
$$
Note that if $\bx$ is a node on the mesh, $\bx_L$ is a node on the lattice.
For derivatives, the subscript $L$ means spatial differentiation with respect to $\bx_L$.
In this setting, Eqs.~\eqref{eq-H} and \eqref{eq-H-source-sink} read
\begin{align}
\partial_{t_L} H_L - \frac{\delta t}{\delta x^2}\Delta_L \theta_L &= 0,  \label{eq:enthalpy_H_LBM} \\
\partial_{t_L} \theta_L - \frac{\delta t}{\delta x^2}\Delta_L \theta_L &= -\frac{1}{\Ste}\partial_{t_L} \varphi(\theta_L).\label{eq:enthalpy_T_LBM}
\end{align}

\subsection{LBM setting}
Let us denote by $\be_i$, $i=0,\cdots,q-1$ the set of discrete velocities for the lattice we want to consider.
We denote by $w_i$ the standard weights associated to the lattice.
We introduce the distribution functions $f_i$, $f_i^{eq}$, $i=0,\cdots,q$; $f_i^{eq}$ is the equilibrium distribution and its value depends on the selected scheme. The lattices associated to a $d-$dimensional problem with $q$ discrete velocities is denoted by D$d$Q$q$. In particular, we will use D1Q3 lattices in 1D and D2Q5 in 2D.

Throughout this note, if $S$ is a function of $t_L$ and $\bx_L$, or if $S_i$ is a distribution function, we set:
\begin{align*}
S &= S(t_L,\bx_L) ,\\
\tilde S &= S(t_L+1,\bx_L) , \\
S_i^* & =S(t_L+1,\bx_L+\be_i), \\
\hat S_i &= S(t_L,\bx_L+\be_i), \\
\check S_i &= S(t_L,\bx_L-\be_i).\\
\end{align*}
We denote by $g_i$ a source term (to be specified later).
Following the Bhatnagar-Gross-Krook (BGK) approximation, the LBM algorithm reads:
\begin{equation}\label{eq:std_lbm_collide_stream}
f_i^* - f_i = -\dfrac{1}{\tau}\left( f_i - f_i^{eq}\right) + g_i,
\end{equation}
where $\tau$ is a relaxation parameter.

For a distribution function $s_i$, we denote by $\calM_k(s)$ the moment of order $k$, i.e.
$$
\calM_k(s) := \sum_{i=0}^{q-1}  \underbrace{\be_i\otimes\cdots\otimes\be_i}_{k\textrm{ times}} s_i.
$$
To alleviate the notations, we also define:
$$
\Pi_k := \calM_k(f),\qquad \Lambda_k:=\calM_k(f^{eq}),\qquad \calG_k:=\calM_k(g).
$$
\subsection{Implicit liquid fraction based method (ILFBM)}
One of the first specific LB models for the total enthalpy equation \eqref{eq-H-source-sink} was introduced in~\cite{jiaung2001lattice}. It relies on the use of a source/sink term. The distribution $f_i$ is defined such that $\theta_L = \calM_0(f)$ and $f_i^{eq} = w_i\theta_L$. The source term $g_i$ is intended to satisfy
$$
g_i = \frac{w_i}{\Ste}\left(\ell - \tilde \ell\right),
$$
where the additional field $\ell$ (resp. $\tilde\ell$) is an approximation of the liquid fraction at time $t$ (resp. $t+\delta t$). The goal is to recover
$
\calG_0 = -\frac{1}{\Ste}\partial_t\varphi(\theta_L),
$
at least, up to third order in $\eps$ (which is the parameter used in the Chapman-Enskog expansion).
$\ell$ is computed from the enthalpy $H$ using:
\begin{equation*}
\ell = l(H) = \begin{cases} 0 &\textrm{ if } H<0, \\ \Ste\ H &\textrm{ if } 0\leq H \leq \frac{1}{\Ste}, \\ 1 &\textrm{ if }H>\frac{1}{\Ste}.\end{cases}
\end{equation*}
The time step procedure to get $\tilde f_i$ and $\tilde \ell$ from $f_i$ and $\ell$ is given in Algorithm~\ref{algo:ilfbm}.
%
The target equation~\eqref{eq:enthalpy_T_LBM} can be recovered from Chapman-Enskog expansion, see \cite{jiaung2001lattice}.
We refer to this approach as the Implicit Liquid Fraction Based Method (ILFBM).

\begin{algorithm}
\caption{Implicit Liquid Fraction Based Method step}\label{algo:ilfbm}
\begin{algorithmic}[1]
\State moments: $\theta\gets\calM_0(f)$ and $H\gets \varphi(\theta)$
\State equilibrium: $f_i^{eq}\gets w_i\theta$
\State inner loop initialization: $\ell^{(0)} \gets \ell$, $f_i^{(0)} \gets f_i$
\While{($\ell^{(k)}$ and $f_i^{(k)}$ not converged)}
    \State collision: $\Omega_i\gets f_i -\frac{1}{\tau}\left( f_i - f_i^{eq}\right) + \frac{w_i}{\Ste}\left(\ell - \ell^{(k-1)}\right)$
    \State streaming: $f_i^{(k)} \gets  \check\Omega_i$
    \State update enthalpy: $H^{(k)} \gets \mathcal M_0(f^{(k)}) + \frac{1}{\Ste}\ell^{(k-1)}$
    \State update liquid fraction : $\ell^{(k)} \gets l(H^{(k)})$
\EndWhile
\State output: $\tilde f_i \gets f_i^{(k)}$ and $\tilde\ell = \ell^{(k)}$
\end{algorithmic}
\end{algorithm}

\subsection{Explicit enthalpy based method (EEBM)}

Another approach was introduced in \cite{huang2013new} to take into account the enthalpy and avoid iterations. $f_i^{eq}$ is modified such that its zeroth order moment is exactly the enthalpy $H$; the other moments are unchanged. $\theta$ can be recovered from $H$ using the following expression:
$$
\theta = \Theta(H) = \begin{cases} %
    H, &\textrm{ if } H < 0, \\ %
    0, &\textrm{ if } 0\leqslant H\leqslant \ds{\frac{1}{\Ste}}, \\ %
    H-\ds{\frac{1}{\Ste}}, &\textrm{ if } H>\ds{\frac{1}{\Ste}}, \end{cases}
$$
hence no source term is required. The time step procedure to get $\tilde f_i$ and $\tilde \ell$ from $f_i$ and $\ell$ is given in Algorithm~\ref{algo:eebm}. No inner loop is required, since $\theta$ can be explicitly computed from $H$.
The target equation~\eqref{eq:enthalpy_H_LBM} can be recovered from Chapman-Enskog expansion, see~\cite{huang2013new}.
We refer to this approach as the Explicit Enthalpy Based Method (EEBM).

\begin{algorithm}
\caption{Explicit Enthalpy Based Method step}\label{algo:eebm}
\begin{algorithmic}[1]
\State moments: $H \gets \mathcal M_0(f)$, $\theta\gets\Theta(H)$
\If {$(i\neq 0)$}
\State equilibrium: $f_i^{eq} = w_i \theta$
\Else
\State equilibrium: $f_0^{eq} = H - (1-w_0)\theta$
\EndIf
\State collision: $ \Omega_i = f_i - \frac{1}{\tau}\left(f_i-f_i^{eq}\right) $
\State streaming: $\tilde f_i \gets  \check\Omega_i$
\end{algorithmic}
\end{algorithm}


\section{Implicit regularized enthalpy based method (IREBM)}\label{sec:Newton_enthalpy}
\subsection{Main idea of the scheme}
In this section, we present our approach to compute the enthalpy, which we will refer to as the Implicit Regularized Enthalpy Based Method (IREBM).
We first replace the discontinuous liquid fraction $\varphi$ in~\eqref{eq:enthalpy_T_LBM} with  the regularized function 
$$
\varphi_\delta(\theta) = \frac12\left(1+\tanh\left(\frac\theta\delta\right)\right),
$$
where $\delta$ is a small positive smoothing parameter (see \cite{sadaka2020parallel}). 
We derive a LB model to recover the regularized equation
\begin{equation}\label{eq:regularized}
\partial_{t_L} \theta_L - \frac{\delta t}{\delta x^2}\Delta_L \theta_L = -\frac{1}{\Ste} \partial_{t_L}\varphi_\delta(\theta_L).
\end{equation}
As for the ILFBM, we will use a source term and we want that $
 \theta_L = \calM_0(f)$ and $f_i^{eq} = w_i \theta_L$.
The idea of the algorithm is to execute an inner loop in order to set
\begin{equation}\label{eq:rhs_fl_1}
g_i = \hat s_i - s_i^*, \qquad\textnormal{with}\qquad s_i = \frac{w_i}{\Ste} \varphi_\delta(\theta_L).
\end{equation}
Note that both explicit ($\hat s_i$) and implicit ($s_i^*$) terms are evaluated at $\bx_L+\be_i$.
The reason is that the use of $s_i^*$ instead of $\tilde s_i$ allows us to recover a local (in space) Newton problem.
Then the use of $\hat s_i$ is required instead of a more usual $s_i$ in order to avoid extra terms in the Chapman-Enskog expansion, see \eqref{eq:CE_exp1}.

\noindent The LBM update scheme is
\begin{equation}\label{eq:irebm_update_scheme}
f_i^* - f_i = -\frac{1}{\tau}\left(f_i - f_i^{eq}\right) + \left( \hat s_i - s_i^*\right).
\end{equation}
The only unknown in the right hand side is $s_i^*$, for which we would need $\tilde\theta_L$.
Let us introduce $\Omega_i$ and $q_i$ such that
\begin{equation}\label{eq:def_qi_omegai}
\hat q_i = \Omega_i:= \left(1-\frac1\tau\right)f_i + \frac1\tau f_i^{eq}.
\end{equation}
Then the previous equation reads
\begin{equation}\label{eq:rhs_fl_stag}
f_i(t_L+1,\bx_L) = q_i(t_L,\bx_L) + \frac{w_i}{\Ste}\left( \varphi_\delta(\theta(t_L,\bx_L)) - \varphi_\delta(\theta(t_L+1,\bx_L))\right).
\end{equation}
Taking the zeroth order moment of \eqref{eq:rhs_fl_stag} yields:
\begin{equation}\label{eq:rhs_fl_2}
\tilde \theta_L = \calM_0\left(q\right) + \frac{1}{\Ste}\left( \varphi_\delta(\theta_L) - \varphi_\delta(\tilde\theta_L)\right).
\end{equation}
We solve this (local in $\bx_L$) non linear equation using a Newton algorithm in order to get $\tilde\theta_L$.
The complete algorithm to get $\tilde f_i$ from $f_i$ is given in Algorithm~\ref{algo:irebm}.

\begin{algorithm}
\caption{Implicit Regularized Enthalpy Based Method step}\label{algo:irebm}
\begin{algorithmic}[1]
\State moments: $\theta\gets \calM_0(f)$
\State equilibrium: $f_i^{eq}\gets w_i \theta$
\State collision: $\Omega_i \gets \left(1-\frac1\tau\right)f_i + \frac1\tau f_i^{eq}$
\State streaming: $q_i \gets \check\Omega_i$
\State Newton initialization: $m_0 \gets \calM_0(q)$, $\theta^{(0)} = \theta$
\While{($\theta^{(k)}$ not converged)}
  \State
    $\ds{\theta^{(k+1)} \gets \frac{m_0 + \frac1{\Ste}\varphi_\delta(\theta) + \frac1{\Ste}\varphi_\delta'(\theta^{(k)})\theta^{(k)}- \frac1{\Ste}\varphi_\delta(\theta^{(k)})}{1+\frac{1}{\Ste}\varphi_\delta'(\theta^{(k)})}}
    $
\EndWhile
\State source term: $g_i \gets \frac{w_i}{\Ste}\left( \varphi_\delta(\theta) - \varphi_\delta(\theta^{(k)})\right)$
\State update: $\tilde f_i = q_i + g_i$
\end{algorithmic}
\end{algorithm}

\subsection{Treatment of boundary conditions}
In order to get a complete LBM scheme in 1D or 2D, one has to impose boundary conditions on the distribution $f_i$. For EEBM and IFLBM, this can be done using standard techniques for Dirichlet or Neumann boundary conditions, cf. \eg \cite{Kaluza2012,ChenMuller2020}. For IREBM, owing to~\eqref{eq:def_qi_omegai} and \eqref{eq:rhs_fl_stag}, we have to set boundary condition on the distribution $q_i$. One can notice that $q_i(t_L,\bx_L)$ is not defined if $\bx_L-\be_i$ is not a point of the lattice, but this value is still required to get the updated distribution function. Since the algorithm depends on the lattice used to solve the numerical problem, we present some examples for simple lattices in 1D and 2D.
Denote by $\Gamma_D$ (resp.
$\Gamma_N$) the part of the boundary on which we want to impose the Dirichlet condition $\theta=\theta_{Dir}$ (resp. the Neumann condition $\partial_n \theta = 0$).

\subsubsection{Dirichlet boundary condition on a D1Q3 or D2Q5 lattice}
Let us consider $\bx_L\in\Gamma_D$ such that $\bx_L$ is not on a corner (the idea can be adapted for the corners). From \eqref{eq:rhs_fl_2} and the definition of the moments, we obtain, for this point $\bx_L$, that
$$
\tilde\theta = \sum_{i=0}^{q-1} q_i + \frac{1}{\Ste}\left(\varphi_\delta(\theta) - \varphi_\delta(\tilde\theta)\right).
$$
Using the Dirichlet boundary condition, this can be rewritten as:
\begin{equation}\label{eq:Dirichlet}
\mathcal M_0(q) = \tilde\theta_{Dir} + \frac{1}{\Ste}\left(\varphi_\delta(\theta_{Dir}) - \varphi_\delta(\tilde\theta_{Dir})\right).
\end{equation}
For a general lattice, one might use any technique to impose a Dirichlet condition on $q_i$, according to \eqref{eq:Dirichlet}. In the case of a D1Q3 or D2Q5 lattice, only one of the $q_i(t,\bx_L)$ is not well-defined in \eqref{eq:def_qi_omegai}.
Let us denote by $i_x$ its index. Then \eqref{eq:Dirichlet} reads 
\begin{equation}\label{eq:Dir_on_q}
q_{i_x} = \tilde\theta_{Dir} + \frac{1}{\Ste}\left(\varphi_\delta(\theta_{Dir}) - \varphi_\delta(\tilde\theta_{Dir})\right) - \sum_{i\neq i_x} q_i.
\end{equation}
Since all the terms in the right hand side are well defined and known, this is a consistent way of defining $q_{i_x}$. 

\subsubsection{Bounce-back condition}

A standard way to impose the homogeneous Neumann boundary condition is to use the so-called "bounce-back" condition \cite{d1986lattice}. For any $\bx_L\in\Gamma_N$, we denote by $i_x$ any index such that $q_{i_x}(t,\bx_L)$ is not well-defined in~\eqref{eq:def_qi_omegai}. The bounce-back condition states that
$$
f_{i_x} = f_{\bar i_x},
$$
at any time, where $\bar i_x$ is the index of the opposite direction to the discrete velocity $e_{i_x}$. Using \eqref{eq:rhs_fl_stag}, this leads to
$$
q_{i_x} + \frac{w_{i_x}}{\Ste}\left( \varphi_\delta(\theta)- \varphi_\delta(\tilde\theta)\right) = q_{\bar i_x} + \frac{w_{\bar i_x}}{\Ste}\left( \varphi_\delta(\theta)- \varphi_\delta(\tilde\theta)\right).
$$
Since $w_{i_x} = w_{\bar i_x}$, this reduces to a bounce-back condition on $q$, i.e. $q_{i_x} = q_{\bar i_x}.$
Since the term in the right-hand side is well defined and known, this is a consistent way of defining $q_{i_x}$. Note also that this bounce-back condition can be used on any lattice.

\subsubsection{Neumann boundary condition on a D2Q5 lattice}

Let us consider $\bx_L\in\Gamma_N$ such that $\bx_L$ is not a corner (the idea can be adapted for the corners). In the particular case of a D2Q5 lattice, only one of the $q_i(t,\bx_L)$ is not well-defined in \eqref{eq:def_qi_omegai}. Let us denote by $i_x$ its index.
An alternative solution to the bounce back condition (with order 1 accuracy in space) is to impose (with our notations) that $\theta(t_L,\bx_L) = \theta(t_L,\bx_L+\be_{i_x})$ for any time.
From \eqref{eq:rhs_fl_2}, we have
\begin{align*}
\tilde\theta(\bx_L) &= \sum_{i=0}^{q-1} q_i(\bx_L) + \frac{1}{\Ste}\left( \varphi_\delta(\theta(\bx_L))-\varphi_\delta(\tilde\theta(\bx_L))\right), \\
\tilde\theta(\bx_L+\be_{i_x}) &= \sum_{i=0}^{q-1} q_i(\bx_L+\be_{i_x}) + \frac{1}{\Ste}\left( \varphi_\delta(\theta(\bx_L+\be_{i_x}))-\varphi_\delta(\tilde\theta(\bx_L+\be_{i_x}))\right).
\end{align*}
The applied Neumann condition at time $t$ ensures that $\theta(\bx_L) = \theta(\bx_L+\be_{i_x})$.
Thus, imposing $\tilde\theta(\bx_L) = \tilde\theta(\bx_L+\be_{i_x})$ is equivalent to satisfy
$$
 \sum_{i=0}^{q-1} q_i(\bx_L) =  \sum_{i=0}^{q-1} q_i(\bx_L+\be_{i_x}).
$$
As a result, we have
\begin{equation}\label{eq:Neu_on_q}
q_{i_x}(\bx_L) = \sum_{i=0}^{q-1}q_i(\bx_L+\be_{i_x}) - \sum_{i\neq i_x} q_i(\bx_L).
\end{equation}
Since all the terms in the right-hand side are well defined and known, this is a consistent way of defining $q_{i_x}(\bx_L)$. 


\section{Numerical analysis of the scheme}\label{sec:analysis}

The goal of this section is to recover the target equation~\eqref{eq:enthalpy_T_LBM} from the LBM setting, using a Chapman-Enskog expansion.
%
We start from the LBM update scheme~\eqref{eq:irebm_update_scheme} and assume that the following holds (standard assumptions for the Chapman-Enskog expansion), for a small parameter $\eps$:
\begin{align*}
f_i &= f_i^{(0)} + \eps f_i^{(1)} + \eps^2 f_i^{(2)}, \\
f_i^{eq} &= f_i^{eq,0}, \\
\partial_{t_L} &= \eps\partial_{t_1} + \eps^2\partial_{t_2}, \\
\nabla_L &= \eps\nabla_{\tilde x}.
\end{align*}
Since we use that $f_i^{eq}=w_i\theta_L$, we obtain 
\begin{align}
\Lambda_0 &= \Pi_0 = \theta_L, \label{eq:lambda0}\\
\Lambda_1 &= 0,\label{eq:lambda1} \\
\Lambda_2 &= \frac13\theta_L\polI.\label{eq:lambda2}
\end{align}  
We also denote by $\Pi_k^{(m)}$ the $k-$th order moment of the distribution $f_i^{(m)}$.
Expanding \eqref{eq:irebm_update_scheme} up to third order in $\eps$, we get:
\begin{align}
\eps D_{1,i} f_i + \frac{\eps^2}{2} D_{1,i}^2 f_i +\partial_{t_2} f_i 
&= -\frac{1}{\tau}\left(f_i-f_i^{eq}\right) - \left( \eps\partial_{t_1}s_i + \eps^2\partial_{t_2}s_i + \frac{\eps^2}{2}\partial_{t_1}^2s_i +  \eps^2\partial_{t_1}\left(\eigrad s_i\right)\right), \label{eq:CE_exp1}
\end{align}
where $D_{1,i}f = \partial_{t_1}f + \eigrad f$.
By identification with respect to consecutive orders of $\eps$, we obtain:
\begin{align}
\eps^0:\ & f_i^{(0)} = f_i^{eq}, \label{eq:CE_fi0}\\
\eps^1:\ & D_{1,i}f_i^{(0)} = -\frac{1}{\tau}f_i^{(1)} - \partial_{t_1}s_i, \label{eq:CE_fi1}\\
\eps^2:\ & \partial_{t_2}f_i^{(0)} + D_{1,i} f_i^{(1)} + \frac{1}{2} D_{1,i}^2 f_i^{(0)} = -\frac{1}{\tau}f_i^{(1)} - \left(\partial_{t_2}s_i + \frac{1}{2}\partial_{t_1}^2s_i +  \partial_{t_1}\left(\eigrad s_i\right)\right).
\label{eq:CE_fi2_orig}
\end{align}
Taking the moments of \eqref{eq:CE_fi0} yields $\Pi_k^{(0)} = \Lambda_k$, for any $k$.
Since we have $\Pi_0=\Lambda_0=\theta_L$, we infer that, for any $l\geq1$:
\begin{equation}\label{eq:pi0k_nul}
\Pi_0^{(l)} = 0.
\end{equation}
Injecting \eqref{eq:CE_fi1} in \eqref{eq:CE_fi2_orig}, we obtain
\begin{equation}\label{eq:CE_fi2}
 \partial_{t_2}f_i^{(0)} + \left(1-\frac1{2\tau}\right)D_{1,i} f_i^{(1)} = -\frac{1}{\tau}f_i^{(1)} - \left(\partial_{t_2}s_i + \frac{1}{2}\partial_{t_1}^2s_i +  \partial_{t_1}\left(\eigrad s_i\right)\right)   + \frac12 D_{1,i}\partial_{t_1} s_i.
\end{equation}

Taking the zeroth and the first order moment of \eqref{eq:CE_fi1} and the zeroth order moment of \eqref{eq:CE_fi2} leads to:
\begin{align}
\partial_{t_1}\theta_L  &= -\frac{1}{\Ste}\partial_{t_1}\varphi_\delta(\theta_L), \label{eq:CE_dt1}\\
\frac13 \divt(\theta_L\polI) &= -\frac{1}{\tau}\Pi_1^{(1)}, \label{eq:CE_dt2}\\
\partial_{t_2}\theta_L + \left(1-\frac1{2\tau}\right)\divt\Pi_1^{(1)} &= -\frac{1}{\Ste}\partial_{t_2}\varphi_\delta(\theta_L). \label{eq:CE_dt3}
\end{align}
Inserting~\eqref{eq:CE_dt2} in~\eqref{eq:CE_dt3} gives
\begin{equation}\label{eq:CE_dt4}
\partial_{t_2}\theta_L - \frac13\left(\tau-\frac12\right)\nabla_{\tilde x}^2\theta_L = -\frac{1}{\Ste}\partial_{t_2}\varphi_\delta(\theta_L).
\end{equation}
We multiply \eqref{eq:CE_dt1} and \eqref{eq:CE_dt4} by $\varepsilon$ and $\varepsilon^2$ respectively and sum these equations to get
\begin{equation}\label{eq:CE_final}
\partial_{t_L}\theta_L - \frac13\left(\tau-\frac12\right)\Delta_L\theta_L = -\frac1{\Ste}\partial_{t_L} \varphi_\delta(\theta_L),
\end{equation}
up to the third order in $\eps$.
We can see that \eqref{eq:CE_final} corresponds to \eqref{eq:enthalpy_T_LBM}, provided that
$$
\tau-\frac12 = \frac{3\delta t}{\delta x^2}.
$$

\section{Numerical illustrations}\label{sec:numerical}

\subsection{1D problem}

\begin{figure}[H]
\centering
\includegraphics[width=0.65\linewidth]{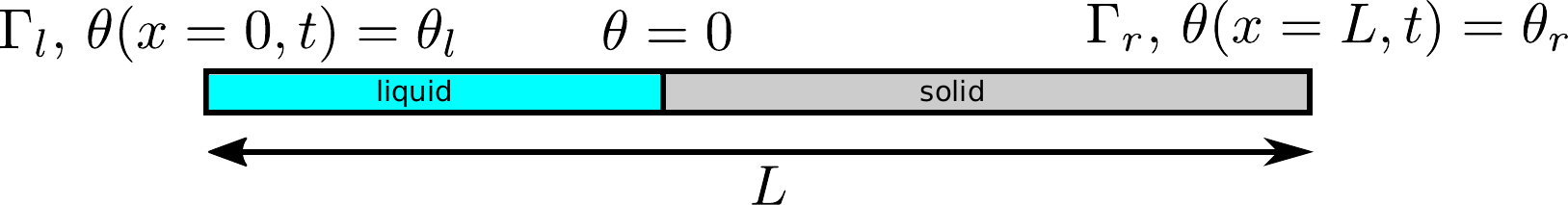}
\caption{Sketch of the simulation configuration.}
\label{fig:case1D_schema}
\end{figure}

To validate the algorithm developed in \S\ref{sec:Newton_enthalpy}, namely the implicit regularized enthalpy based method (IREBM), we consider a case from \cite{huang2013new}.
We compare results obtained by two other algorithms, the explicit enthalpy based method (EEBM) from \cite{huang2013new} and the Implicit Liquid Fraction Based Method (ILFBM) from \cite{jiaung2001lattice}. We use the numerical methods described in \S\ref{sec:problem} and the exact solution as reference.
We consider a one dimensional bar of length $L$ made of aluminum.
The temperature is fixed at both ends, with $T_h$ on the left and $T_c$ on the right.
During the simulation, physical quantities are rendered dimensionless, and dimensionless temperature is denoted by $\theta$.
The case is initialized with a temperature distribution satisfying the Stefan condition for the case where most of the bar is solid. More precisely, the position of the interface is given by $x_f(t)=2\lambda\sqrt{t}$, where $\lambda$ is computed from~\eqref{eq-Stefan-lambda}. We choose $t_0$ such that $x_f(t=t_0) = 0.01$. Starting from the initial temperature distribution given by~\eqref{eq-exact-sol}, we track the displacement of the solid-liquid interface.
One can refer to Figure~\ref{fig:case1D_schema} for a sketch of the simulation; physical properties are given in Table~\ref{tab:case1D_physical_quantities}.

\begin{figure}[H]
\centering
\includegraphics[width=0.9\linewidth]{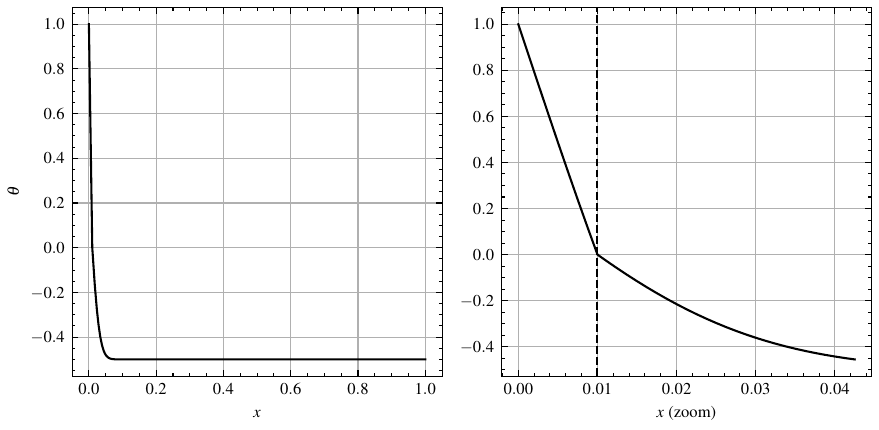}
\caption{Initial solution for the 1D case.}
\label{fig:case1D_initial_solution}
\end{figure}

\begin{table}[H]
    \centering
    \begin{tabular}{|l|l|l|l|}
    \hline
    physical quantity & symbol & value & unit \\
    \hline
    density & $\rho_l=\rho_s$ & 2500 & $\text{kg}\cdot\text{m}^{-3}$ \\
    specific heat capacity & $c_l=c_s$ & 1 & $\text{kJ}\cdot\text{kg}^{-1}\cdot\text{K}^{-1}$ \\
    thermal conductivity & $k_l=k_s$ & 0.2 & $\text{kW}\cdot \text{m}^{-1}\cdot\text{K}^{-1}$\\
    thermal diffusivity & $\alpha_l=\alpha_s$ & $8\times10^{-5}$ & $\text{m}^2\cdot\text{s}^{-1}$\\
    latent heat & $h_{sl}$ & 350 & $\text{kJ}\cdot\text{kg}^{-1}\cdot\text{K}^{-1}$ \\
    fusion temperature & $T_f$ & 500 & $\degree\text{C}$ \\
    temperature left side & $T_h$ & 600 & $\degree\text{C}$ \\
    temperature right side & $T_c$ & 450  & $\degree\text{C}$ \\
    reference length & $L$ & 1 & m \\
    \hline
    \end{tabular}
    \caption{Physical quantities and values used for the 1D case.}
    \label{tab:case1D_physical_quantities}
\end{table}
The non-dimensional formulation of this problem is obtained by using the as reference temperature the difference $\delta T = T_h-T_f$, resulting in the Stefan number $\Ste = 0.2857$. 
%
%
%
%
%
The initial temperature profile is shown in Fig.~\ref{fig:case1D_initial_solution} and exhibits a discontinuity of derivatives at $x_f$, as expected since the phase change occurs at this location.

\begin{figure}[H]
\centering
\includegraphics[width=0.95\linewidth]{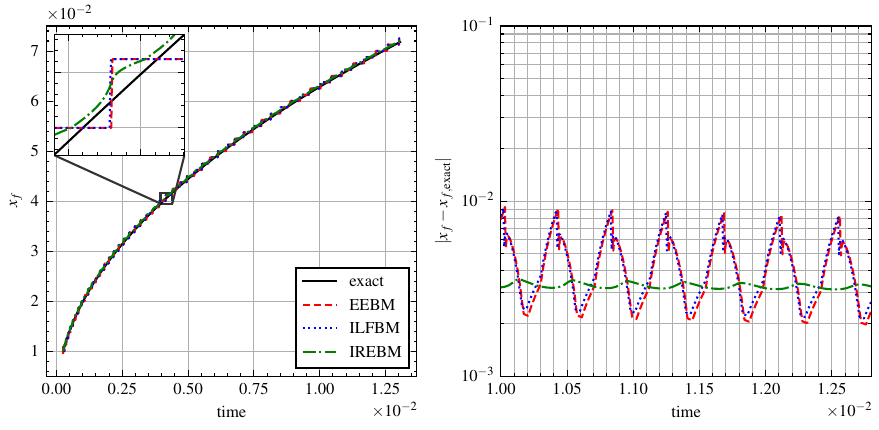}
\caption{Position of the interface (left) and error with respect to the exact solution (right). Comparison between EEBM, ILFBM, IREBM (with $\delta=0.005$) $N = 801$.}
\label{fig:case1D_comp_method_801}
\end{figure}

In Fig.~\ref{fig:case1D_comp_method_801}, we compare the three methods using the same spatial and time resolutions.
EEBM and ILFBM behave almost identicaly, resulting in an oscillation of the error on the position of the interface between $2\times 10^{-3}$ and $1\times 10^{-2}$.
For the implicit regularized enthalpy based method (IREBM), the error on the position of the interface oscillates between $3\times 10^{-3}$ and $4\times 10^{-3}$.
Note that the error fluctuates a lot less with IREBM, which indeed is the result of the added regularization.
One needs to keep in mind that the problem is modified when using the regularization, and methods EEBM, ILFBM will not converge to the solution of IREBM if $\delta$ is kept constant.

\begin{figure}[H]
\centering
\includegraphics[width=0.95\linewidth]{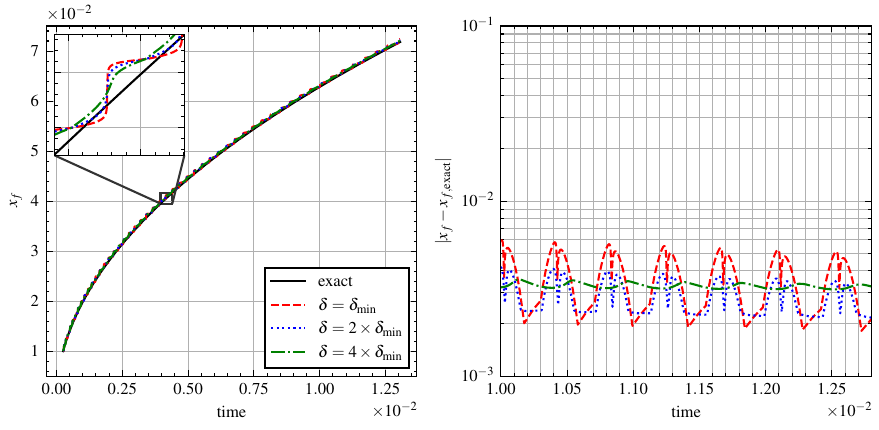}
\caption{Position of the interface (left) and error with respect to the exact solution (right). Comparison of IREBM for three values of $\delta$ ($\delta_\text{min}=0.005$), for $N = 801$.}
\label{fig:case1D_comp_eps_801}
\end{figure}

In Fig.~\ref{fig:case1D_comp_eps_801}, we can see the effect of $\delta$ on the evolution of the position and the oscillations of error of the interface.
Larger is $\delta$, smaller the oscillations of the error are.
If $\delta$ is too large, the regularization will offset the solution to another problem, so $\delta$ needs to be tuned in order to keep the solution close to the original problem.

\begin{figure}[H]
\centering
\includegraphics[width=0.95\linewidth]{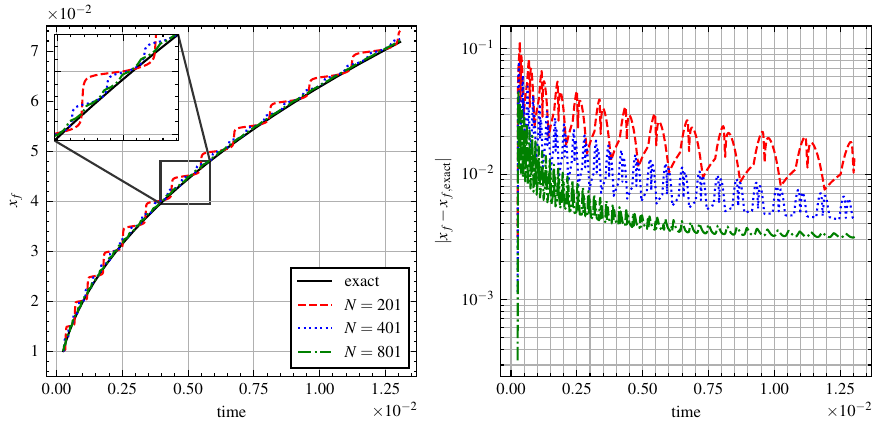}
\caption{Position of the interface (left) and error with respect to the exact solution (right). Comparison of different mesh resolutions for IREBM, with $\delta=0.005$.}
\label{fig:case1D_comp_res_eps5e-3}
\end{figure}

In Fig.~\ref{fig:case1D_comp_res_eps5e-3}, we can see the effect of the resolution on the solution, when keeping $\delta$ constant.
We note that the error on the position of the interface diminishes as the resolution $N$ increases, and the oscillations are reduced.

\begin{figure}[H]
\centering
\includegraphics[width=0.95\linewidth]{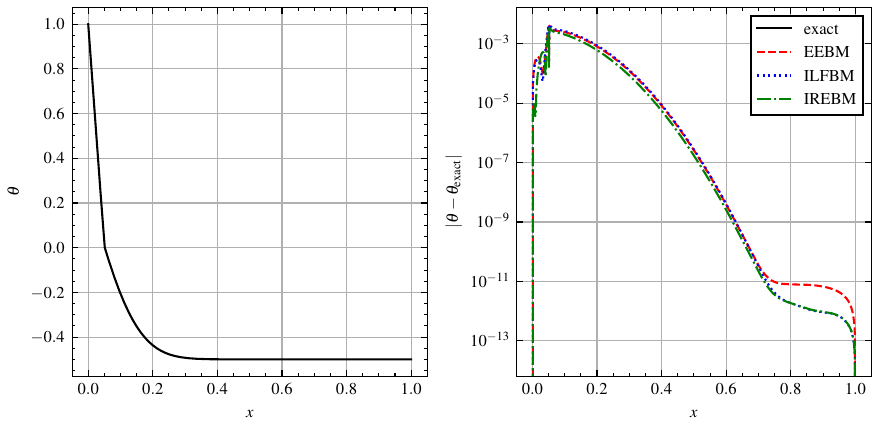}
\caption{Solution (left) and error (right) at final time ($t=160s$) for all the methods (EEBM, ILFBM, IREBM with $\delta=0.005$) and $N = 801$.}
\label{fig:case1D_visu_solution}
\end{figure}

In Fig.~\ref{fig:case1D_visu_solution}, we compare the three methods with the same spatial and time resolution.
All methods exhibit very similar behaviors.
As expected, the maximum of the error is near the phase change interface.
The IREBM is having a slightly lower error than the two other methods.
The maximum of the error is $4.06\times10^{-3}$for the EEBM, $4.50\times10^{-3}$ for the ILFBM and $3.41\times10^{-3}$ for the IREBM.

Finally, we report the impact of Newton iterations on the computational time in Tab.~\ref{tab:computational-time}. As expected, IREBM is slower than EEBM because of Newton iterations. However, IREBM is significantly faster than ILBM for two reasons: (i) the Newton algorithm usually requires a low number of iterations to converge; (ii) inner iterations are performed on only one field, while the ILFBM inner loop involves all $q$ fields of the particle distribution function and a streaming process. For these reasons, the extra cost of the Newton loop is expected to be even lower in 2D or 3D, when larger values of $q$ in the lattice are necessary.

\begin{table}[h]
\centering
\begin{tabular}{|c||c|c|c|}\hline
 $N$ & EEBM & IREBM & ILFBM \\ \hline
 201 & 0.026 &  0.106 & 0.422 \\
 401 & 0.200 & 0.739 & 2.997\\
 801 & 1.538 &  5.586 & 22.01\\
 1601 & 12.43 & 43.79 & 164.06 \\ \hline
\end{tabular}
\caption{1D case, lattice D1Q3: computational time (in seconds, sequential code on Mac M1), with $\tau=0.62$ and $\delta=0.01$ for IREBM.}%
\label{tab:computational-time}%
\end{table}

\subsection{2D problem}

\begin{figure}[H]
\centering
\includegraphics[width=0.25\linewidth]{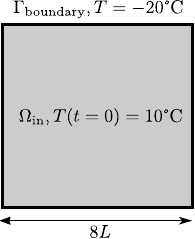}
\caption{Sketch of the 2D simulation configuration.}
\label{fig:case2D_schema}
\end{figure}

\begin{figure}[H]
\centering
\includegraphics[width=0.95\linewidth]{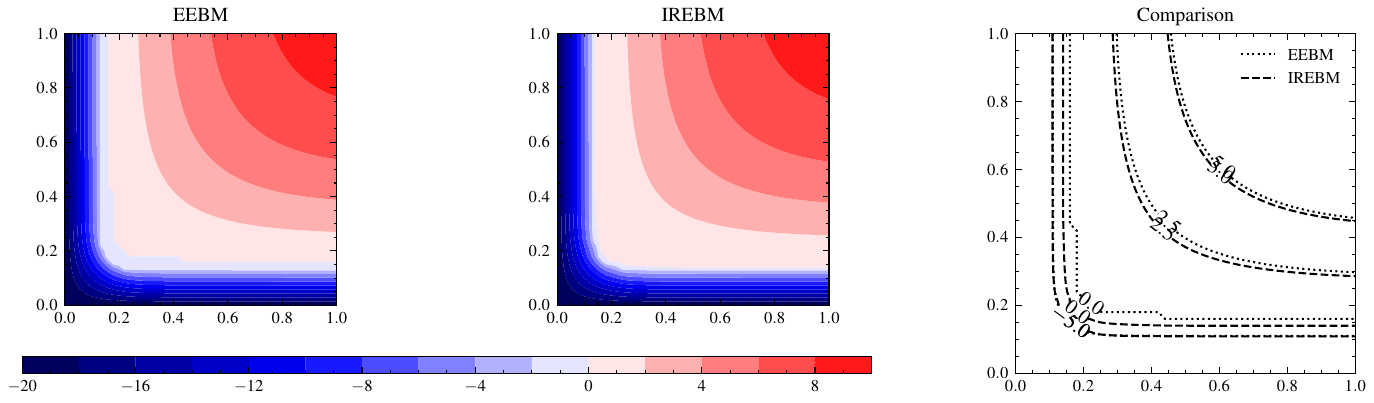}
\caption{2D case: distribution of temperature for EEBM (left) and IREBM (center) and position of three isovalues (right). Simulations for $N=201$ and $\delta=2\times 10^{-2} \approx \Delta x$.}
\label{fig:case2D_comparison_iso_201}
\end{figure}

In order to show the 2D capabilities of the developed algorithm, we compare our results to \cite{prapainop2004simulation}.
They based their comparison on an approximate \textit{analytical} solution from \cite{jiji1970two}.
Other analysis of solutions at corners is given in \cite{king1999two}.
The setup consists of a squared cavity of $8m$ wide, initially at $T = T_i =10\degree\text{C}$ and cooled at all sides to $T = T_\Gamma =-20\degree\text{C}$.
Using the symmetries of the problem, we solve Eq.~\ref{eq-H-source-sink} in a square representing the lower-left quarter of the original domain, see Fig.~\ref{fig:case2D_comparison_iso_201}.
One can refer to Fig.~\ref{fig:case2D_schema} for a sketch of the simulation; physical properties are given in Tab.~\ref{tab:case2D_physical_quantities}.

The computational domain (see e.g. Figure~\ref{fig:case2D_comparison_iso_201}) is then a square, in which Eq.~\eqref{eq-H-source-sink} is solved. The bottom and left walls are kept at the cold temperature, resulting in a Dirichlet boundary condition. To take into account symmetry properties, the system is supplemented with the homogeneous Neumann condition on the top and right walls.

\begin{table}[H]
    \centering
    \begin{tabular}{|l|l|l|l|}
    \hline
    physical quantity & symbol & value & unit \\
    \hline
    density & $\rho_l=\rho_s$ & 1000 & $\text{kg}\cdot\text{m}^{-3}$ \\
    specific heat capacity & $c_l=c_s$ & 1.762 & $\text{kJ}\cdot\text{kg}^{-1}\cdot\text{K}^{-1}$ \\
    thermal conductivity & $k_l=k_s$ & 2.22$\times10^{-3}$ & $\text{kW}\cdot \text{m}^{-1}\cdot\text{K}^{-1}$\\
    thermal diffusivity & $\alpha_l=\alpha_s$ & $1.26\times10^{-6}$ & $\text{m}^2\cdot\text{s}^{-1}$\\
    latent heat & $h_{sl}$ & 338 & $\text{kJ}\cdot\text{kg}^{-1}\cdot\text{K}^{-1}$ \\
    fusion temperature & $T_f$ & 0 & $\degree\text{C}$ \\
    initial temperature & $T_i$ & 10 & $\degree\text{C}$ \\
    temperature boundary & $T_\Gamma$ & -20  & $\degree\text{C}$ \\
    reference length & $L$ & 1 & m \\
    \hline
    \end{tabular}
    \caption{Physical quantities and their values used for the 2D case.}
    \label{tab:case2D_physical_quantities}
\end{table}

The non-dimensional formulation of this problem is obtained by using as reference temperature the difference $\delta T = T_i-T_f$, resulting in the Stefan number $\Ste = 0.0521$. 
Figure~\ref{fig:case2D_comparison_iso_201} compares both EEBM and IREBM (ILFBM did not run successfully) at a final time $t=20$h.
Both methods behave similarly, excepting at the interface where EEBM, as expected, has an oscillatory behavior. The level set $T=-5\degree\text{C}$ is identical for both methods, but positive isovalues are more distinct.

\begin{figure}[H]
\centering
\includegraphics[width=0.95\linewidth]{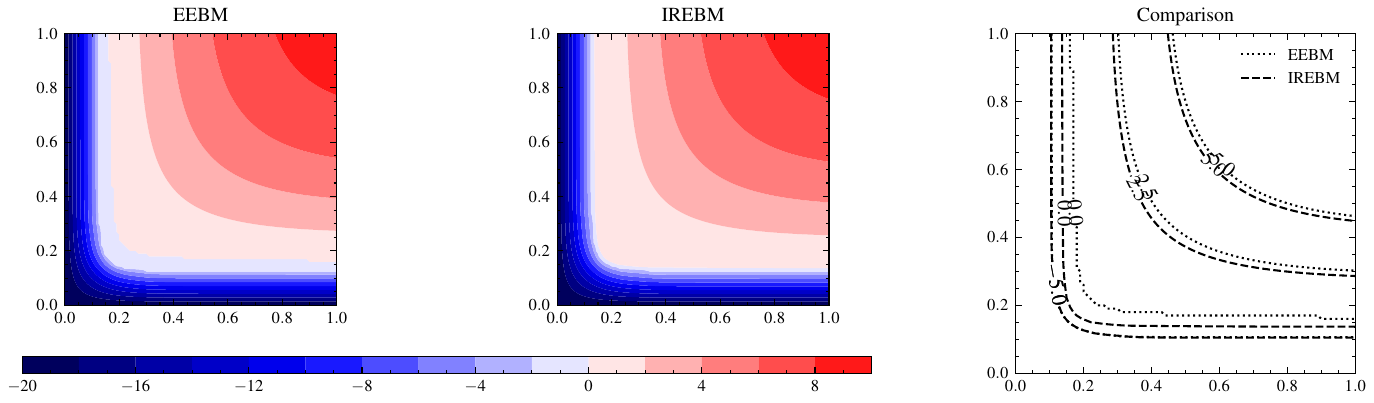}
\caption{2D case: distribution of temperature for EEBM (left) and IREBM (center) and position of some isovalues (right). Simulations for $N=401$ and $\delta=2\times 10^{-2}\approx 2 \Delta x$.}
\label{fig:case2D_comparison_iso_401}
\end{figure}

To assess that the displacement of the isovalues is not due to the spatial resolution, Fig.~\ref{fig:case2D_comparison_iso_401} compares both EEBM and IREBM at a final time $t=20$h but for a higher resolution than previously ($n_{x,y}=401$ compared to $n_{x,y}=201$).
The regularization value $\delta$ is kept the same ($\delta=2\times 10^{-2}$).
The same overall behavior is observed for both methods, with the same oscillating behavior for EEBM and $T=0\degree\text{C}$.

\begin{figure}[H]
\centering
\includegraphics[width=0.55\linewidth]{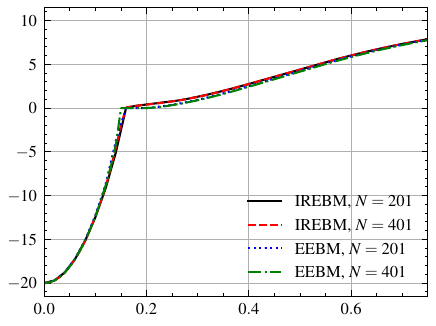}
\caption{2D case: temperature along the $x=y$ axis at $t=20h$. Comparison of EEBM and IREBM for different mesh resolutions.}
\label{fig:case2D_xy_comparison}
\end{figure}

\begin{figure}[H]
\centering
\includegraphics[width=0.55\linewidth]{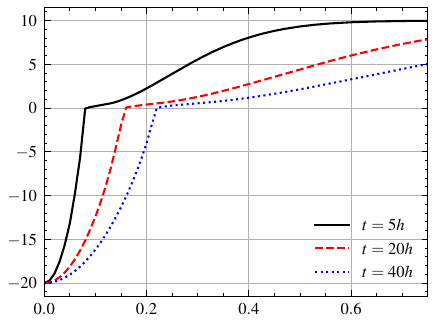}
\caption{2D case: temperature along the $x=y$ axis for $t=5$, $20$, $40$h (IREBM simulation, $N=401$).}
\label{fig:case2D_xy_several_t}
\end{figure}

Figures~\ref{fig:case2D_xy_comparison} and~\ref{fig:case2D_xy_several_t} show the temperature along the $x=y$ axis. In Fig.~\ref{fig:case2D_xy_comparison}, EEBM and IREBM are compared for two discretizations ($N=201$ and $N=401$), and behave similarly. In Fig.~\ref{fig:case2D_xy_several_t}, several time steps are shown. Our results match those of \cite{prapainop2004simulation}.

Finally, we report the impact of the Newton iterations on the computational time in Tab.~\ref{tab:computational-time-2D}. The IREBM is slower than the EEBM with a close to 3 for the D2Q5 lattice and 2.5 for the D2Q9 lattice. This is indeed better than the 1D case (see Tab.~\ref{tab:computational-time}) where the ratio was close to 3.5. As expected, for large values of $q$, the extra cost for the Newton iterations is reduced. This is explained by the fact that the Newton iterations are performed on only one field, independently of the value of $q$. Note that ILFBM results are not shown here because the computational time is significantly larger. The reason is that, for all spatial resolutions and time steps considered here, the inner iterations failed to converge, and the scheme does not yield a proper solution.

Note also that, from a strict computational point of view, when $N$ is multiplied by 2, the computational cost of a LBM scheme for a 2D case should be multiplied by 16. Here we notice that going from $N = 201$ to $N = 401$, the ratio is greater than 25. This is due to the fact that for $N=101$ and $N=201$, all the lattice stands in the L2 memory, but for $N=401$ this is no longer the case (12MB on Apple M1). 

\begin{table}[h]
\centering
\begin{tabular}{|c||c|c|c||c|c|c|}\hline
 & \multicolumn{3}{c||}{D2Q5} & \multicolumn{3}{c|}{D2Q9} \\ \hline
 $N$ & EEBM & IREBM & ratio & EEBM & IREBM & ratio \\ \hline
 101 & 0.064 & 0.250 & 3.90 & 0.104 & 0.303 & 2.91 \\
 201 & 0.875 & 3.794 & 4.34 & 1.697 & 5.133 & 3.02 \\
 401 & 25.194 & 74.245 & 2.95 & 43.452 & 106.58 & 2.45 \\
 801 & 656.58 & 1675.79 & 2.55 & 1028.56 & 2188.74 & 2.13 \\ \hline
\end{tabular}
\caption{2D case: computational time (in seconds, sequential code on Mac M1), with $\tau=0.84$ and $\delta=0.01$ for IREBM.}%
\label{tab:computational-time-2D}%
\end{table}

\section{Conclusion}

A new implicit regularized enthalpy based scheme was developed to solve the Stefan problem using the LBM formalism. The new method was analysed using both Taylor expansion and Chapman-Enskog analysis.
This new approach (IREBM) incorporates ideas from both the implicit liquid fraction based method and also from the explicit enthalpy based method.
The regularization introduced in the method reduces the oscillations of the position of the interface observed in the results obtained with previously published methods.
Considering the position of the interface, the results obtained with IREBM show a better stability and lower error when compared to the two other methods.
At a given time, the overall solution is similar to that obtained with EEBM or ILFBM.
Parallelization and vectorization is kept identical and the arithmetic intensity is augmented due to the Newton loop.
Future work will include new capabilities, new test cases on 3D configurations, and GPU open-source implementation.

\section*{Acknowledgements}

This project was co-financed by the European Union with the European regional development funds and by the Normandy Regional Council via the ELBA project. Part of this work was performed using computing resources of CRIANN (Centre R\'egional Informatique et d'Applications Num\'eriques de Normandie, France). The authors also thank Pr. U. R\"ude and his group at FAU (Friedrich-Alexander-Universit\"at) for fruitful discussions and insights regarding the detailed implementation of LBM.

\bibliographystyle{elsarticle-num}
\bibliography{stefan}

\end{document}